\pdfoutput=1
\RequirePackage{ifpdf}
\ifpdf 
\documentclass[pdftex]{sigma}
\else
\documentclass{sigma}
\fi

\numberwithin{equation}{section}

\usepackage{cite}

\begin{document}

\allowdisplaybreaks

\renewcommand{\thefootnote}{$\star$}

\renewcommand{\PaperNumber}{016}

\FirstPageHeading

\ShortArticleName{A Generalization of the Hopf--Cole Transformation}

\ArticleName{A Generalization of the Hopf--Cole Transformation\footnote{This
paper is a~contribution to the Special Issue ``Geometrical Methods in Mathematical Physics''.
The full collection is available
at \href{http://www.emis.de/journals/SIGMA/GMMP2012.html}{http://www.emis.de/journals/SIGMA/GMMP2012.html}}}

\Author{Paulius MI\v{S}KINIS}

\AuthorNameForHeading{P.~Mi\v{s}kinis}

\Address{Department of Physics, Faculty of
Fundamental Sciences,
\\
Vilnius Gediminas Technical University, Saul\.{e}tekio Ave 11, LT-10223, Vilnius-40, Lithuania}
\Email{\href{mailto:paulius.miskinis@vgtu.lt}{paulius.miskinis@vgtu.lt}}

\ArticleDates{Received June 04, 2012, in f\/inal form February 17, 2013; Published online February 25, 2013}

\Abstract{A generalization of the Hopf--Cole transformation and~its relation
to the Bur\-gers equation of integer order and~the dif\/fusion
equation with quadratic nonlinearity are discussed.
The explicit form of a~particular analytical solution is presented.
The existence of the travelling wave solution and~the interaction of
nonlocal perturbation are considered.
The nonlocal generalizations of
the one-dimensional dif\/fusion equation with quadratic nonlinearity and
of the Burgers equation are analyzed.}


\Classification{26A33; 35K55; 45K05}

\renewcommand{\thefootnote}{\arabic{footnote}}
\setcounter{footnote}{0}

\section{Introduction}

The classical Hopf--Cole transformation is applied
for the solution of the nonlinear dif\/fusion equation.
There are
two well-known nonlinear generalizations of the dif\/fusion
equation: with quadratic nonlinearity and~the Burgers equation.
The f\/irst of them has applications in dealing with plasma and~acoustic
phenomena~\cite{GalaktionovPosashkov1989}.
The Burgers equation
was initially proposed by H.~Bateman while modelling the weak
viscous liquid motion~\cite{Bateman1915} and~later rediscovered
by Burgers as a~simple nonlinear partial dif\/ferential equation in
studies on turbulence~\cite{Burgers1974}.
This equation can be
viewed as a~simplif\/ied version of the Navier--Stokes equation and
related to the heat equation via the Hopf--Cole transformation~\cite{Hopf1950,Cole1951}.
Presently, the number of applications of
the Burgers equation is immense (see, for
instance,~\cite{GurbatovMalakhovSaichev1991,Woyczynski1998} and~references
below).

In the case when the properties of a~system in a~certain point of
conf\/iguration or phase space depend not only on the properties
of this system at this point, but also on the properties of at
least one point of the environment, we deal with the nonlocal
phenomena.
As the examples, let us mention the well-known
prey--predator system of Volterra with delay in ecology~\cite{Volterra}, the ferromagnetic
properties of matter in physics
\cite{Bertotti}, viscoelastic phenomena in mechanics~\cite{Lakes}.
From the mathematical point of view, such phenomena are usually
described by the integro-dif\/ferential equations~\cite{AgarwalORegan2001}.
Over the last few years, more attention
has been given to a~special part of the theory of
integro-dif\/ferential equations, the so-called fractional calculus
\cite{OldhamSpanier1974,SamkoKilbasMarichev1993,Podlubny1999}.
This approach is applied not only in the theory of fractals and~to
the above-mentioned, already classical, nonlocal phenomena, but
also for the description of electrical, biological and~dif\/fusion
phenomena.
The latter topic, as follows from the growing number of
publications, receives much attention~\cite{Podlubny1999}.

Some time ago, two fractional generalizations of the classical
dif\/fusion equation were proposed.
One of them leads to replacing the second space derivative by the fractional one
\begin{gather}
\label{fbe:1}
\phi_t-\alpha\,_{a}D^{2+p}_x\phi=0,
\end{gather}
where $\phi=\phi(x,t)$, $_{a}D^p_x$ is a~fractional derivative
in the sense of Caputo~\cite{SamkoKilbasMarichev1993,Podlubny1999},
where $0<p<1$, and~$a$ is a~parameter of nonlocality.

From the physical point of view, we may consider this spatial
fractional derivative as a~Fourier transformation of the
fractional power of the wave number~$k$ (see a~short review in
Appendix~\ref{appendixB}).

From this approach, but for the wave equation, the fractional
derivative was considered by A.N.~Gerasimov in~\cite{Gerasimov}.
The other generalization, proposed by R.R.~Nigmatullin in~\cite{Nigmatullin1984}, is related to
the fractional substitution
of the time derivative
\begin{gather*}
\,_{a}D^p_t\phi-\alpha\phi_{xx}=0.
\end{gather*}

In this situation, due to the well-known relation between the
dif\/fusion and~the Burgers equations, we may expect two analogous
nonlocal generalizations of the Burgers equation
\begin{gather}
\label{fbe:4}
\phi_t+\phi\phi_x-\alpha\,_{a}D^p_x\phi=0,
\end{gather}
and the other one
\begin{gather}
\label{fbe:5}
\,_{a}D^p_t\phi+\phi\phi_x-\alpha\phi_{xx}=0.
\end{gather}
These generalizations~\eqref{fbe:4} and~\eqref{fbe:5} are just
the applications of the ideas of~\cite{Gerasimov,Nigmatullin1984}.
Usually, equation~\eqref{fbe:5} is called the fractional
Burgers equation.
However, most interesting and~perhaps most
productive is the third, fractional Burgers equation with nonlocal
nonlinearity (FBENN), i.e.\ the nonlinear and~nonlocal
generalization of the dif\/fusion equation, based on a~fractional
generalization of the Hopf and~Cole transformation
\begin{gather}
\label{fbe:0}
\phi_t+\frac{1}{2}\,_{a}D^p_x\big({}_{a}D^{1-p}_x\phi
\big)^2-\alpha\phi_{xx}=0,
\end{gather}
in which $\phi(x,t),\phi_0(x)\in\mathbb{R}$,
$-\infty<x<+\infty$; $t\ge0$ and~the parameter $\alpha>0$.
In
this equation,
$_{a}D^{p}_x\phi=\lambda^{p-1}\,_{a}\partial^{p}_x\phi$ is the fractional derivative,
where $\lambda$ is the length parameter, and
$_{a}\partial^{p}_x\phi$ is the Caputo fractional derivative
(see Appendix~\ref{appendixB}).
Thus, all the terms of the FBENN have integer physical
dimensions.
From the mathematical point of view, in the case of
linear equation, the transition to the dimensionless form is
recommended, and~in the case of nonlinear equation such transition
is necessary (for details, see Appendix~\ref{appendixA}).

A series of exact analytical solutions of this equation, the
asymptotic form of the solutions and~a~fractional generalization
of the Reynolds number are presented.
Concrete examples corresponding
to the simplest behavior of fractal solution are
analyzed.

\section{Linearization}
\label{section2}

The main problem in solving the FBENN~\eqref{fbe:0} is its
nonlinearity.
If for the linear integro-dif\/ferential equation
we may apply some powerful methods~\cite{Polyanin_Manzhirov1998,Sakhnovich_1996}, in the nonlinear
case we can get only general estimates~\cite{Kilbas2003}.
Therefore, the application of the old method of linearization
could be helpful in understanding some properties of the FBENN.

Let $\bar{\phi}$ be a~known solution of the FBENN~\eqref{fbe:0}, and
$\phi=\bar{\phi}+\varepsilon\psi$ is a~small perturbation of
solution $\bar{\phi}$, where $\psi=\psi(x,t)$.
Then $\psi(x,t)$, in
the f\/irst order of $\varepsilon$, obeys the linear equation
\begin{gather}
\label{FBE:21}
\psi_t+\bar{\phi}_x{}\,_{a}D^{1-p}_x\psi-\alpha
\psi_{xx}=0.
\end{gather}

The solution $\bar{\phi}$ could be obtained by any method
presented in the next sections.
Thus, according to the general theorems
of the uniqueness and~existence of the weak sense solution
given in equation~\eqref{FBE:05} below, the solution of the
linearized FBENN~\eqref{FBE:21} could be found numerically.

The asymptotic behavior of the linearized solution $\psi(x,t)$ at
$t\rightarrow\pm\infty$ is very important for the stability of
solution $\phi(x,t)$.

\section{Nonlocal perturbations of the local solutions}
\label{section3}

Let us take advantage of the fact that in the case of $p=0$ and
$p=1$ the solutions of the respective equations are known.
To f\/ind
how the nonlocality changes of the local solution of the FBENN
\eqref{fbe:0}, we shall perform the expansion of the latter two
equations in the neighbourhood of $0+\varepsilon$ and
$1-\varepsilon$ of the parameter $p$.

In the neighbourhood of the point $p=0+\varepsilon$, the
fractional derivative has the expansion
\begin{gather*}
\,_{a}D^{p}_x\psi(x)=\phi(x)+\varepsilon\,_{a}\hat{N}_x\phi(x)+O\big(p^2\big),
\end{gather*}
here $_{a}\hat{N}_x\phi(x)$ is the nonlocal operator
\begin{gather*}
\,_{a}\hat{N}_x\phi(x)=- \gamma\phi(x)-\int^x_a\phi '(\xi)\log{(x-\xi)}\,d\xi,
\end{gather*}
here $\gamma$ is the Euler constant.

Then the expansion of the FBENN~\eqref{fbe:0} in the neighbourhood of
the respective points is
\begin{gather}
\label{FBE:33}
\phi_t+\frac{1}{2}\phi^2_x-\alpha
\phi_{xx}=\varepsilon\left[\phi_x \big({}_{a}\hat{N}_x\phi \big)_x-
\frac{1}{2}\,{}_{a}\hat{N}_x\phi^2_x\right],\qquad
p=0+\varepsilon,
\\
\label{FBE:34}
\phi_t+\phi\phi_x-\alpha\phi_{xx}=(1-\varepsilon)\big[{}_{a}\hat{N}_x\phi^2-
\phi\,{}_{a}\hat{N}_x\phi\big]_x,\qquad p=1-\varepsilon.
\end{gather}

\looseness=-1
Thus, as follows from equations~\eqref{FBE:33} and~\eqref{FBE:34},
in the case of weak nonlocality (small values of the parameters
$p$ or $1-p$) FBENN~\eqref{fbe:0} can be interpreted as the
respective classical equations perturbed by nonlocal terms of a
highly specif\/ic form.
Both these equations can be solved
numerically.

\section{The travelling wave solution}
\label{section4}

Let us change the reference frame and~turn to a~new variable $\xi
=x-ut$.
Then the FBENN~\eqref{fbe:0} takes the form
\begin{gather}
\label{FBE:11}
\frac{1}{2}\,{}_{a}D_\xi^p\big({}_{a}D_\xi^{1-p}\phi\big)^2=\alpha\phi''+u\phi'.
\end{gather}

The solutions of the corresponding equations for $p\to0$
and $p\to1$ are known
\begin{gather}
\phi (\xi )=2u\xi-2\alpha\ln\big[{e^{\frac{u\xi}{\alpha}+2uC_{1}}-1}\big]+C_2,
\nonumber
\\
\label{FBE:13}
\phi (\xi )=u\left(1-\textrm{th}\left[{\frac{u}{2\alpha} ({\xi+C} )}
\right]\right).
\end{gather}

Suppose that according to the ``intermediate'' character of the
evolution equation~\eqref{FBE:11}, the travelling wave solution
$\phi(\xi)$ of this equation is just a~fractional
 ``deformation'' of the solution $\phi^{(0)}$ $(\xi
) ({\phi^{(1)}(\xi)}
)$ with the integer value of the parameter $p=0$  $(
 p=1 )$
\begin{gather}
\label{FBE:14}
\phi (\xi )={}_{a}D_\xi^p\phi^{(0)}(\xi),\qquad\phi(\xi
)={}_{a}D_\xi^{p-1}\phi^{(1)}(\xi),
\end{gather}
where $_{a}D_\xi^{p-1}f(\xi)\equiv {}_{a}I_\xi
^{1-p}f(\xi)$ means a~fractional integral (see Appendix~\ref{appendixB}).
Indeed, the substitution of expressions~\eqref{FBE:14} into travelling wave
equation~\eqref{FBE:11} leads to equations
for the corresponding
solutions $\phi^{(0)}(\xi)$ and~$\phi^{(1)}(\xi)$.

As a~consequence, from~\eqref{FBE:14} we get $\phi^{(1
)}(\xi)=\partial_\xi\phi^{(0
)}(\xi)$.
The assumption~\eqref{FBE:14} allows
us to get the travelling wave solution of the equation
\eqref{FBE:11} by substituting the corresponding solutions of this
equation for the limit cases $p\rightarrow0$ and~$p\rightarrow
1$ into expression~\eqref{FBE:14}.

The travelling wave solution of~\eqref{FBE:11} is
\begin{gather*}
\phi=2u\xi-2\alpha\log{\left[\exp{\frac{u(c+\xi)}{\alpha}}-1\right]},\qquad
p=0,
\end{gather*}
with the asymptotics $\phi=-2uc$, $\xi\gg1$;
\begin{gather*}
\phi=-\frac{2}{\Gamma(2-p)}\frac{u\xi+\frac{p}{\alpha}u^2-(1-p)C}{(\xi-a)^p},\qquad
0<p<1,
\end{gather*}
with the integration constant $C$, and
\begin{gather*}
\phi=\phi_1+\frac{\phi_2-\phi_1}{1+\exp{\left(\frac{\phi_2-\phi_1}{2\alpha}\xi\right)}},\qquad
p=1,
\end{gather*}
with the asymptotics $\phi(\xi\rightarrow+\infty)=\phi_1$,
$\phi(\xi\rightarrow-\infty)=\phi_2$,
$\phi_2>\phi_1$.

\section[The fractional generalization  of the Hopf-Cole transformation]{The fractional generalization  of the Hopf--Cole transformation}
\label{section5}

In the case of the evolutional Burgers equation (BE),
\begin{gather}
\label{FBE:04A}
\phi_{t}+\phi\phi_{x}-\alpha\phi_{xx}=0,
\end{gather}
well known is the simple nonlinear Hopf--Cole transformation
\begin{gather}
\label{FBE:61}
\phi(x,t)=-2\alpha\frac{w_x(x,t)}{w(x,t)},
\end{gather}
which relates any solution $w(x,t)$ of the dif\/fusion equation (DE)
$w_t=\alpha w_{xx}$ to the solu\-tion~$\phi(x,t)$ of the BE.

In the case of the nonlinear dif\/fusion equation (NDE) with
quadratic nonlinearity
\begin{gather}
\label{FBE:03B}
\phi_{t}+\frac{1}{2} \phi^{2}_{x}-\alpha
\phi_{xx}=0,
\end{gather}
there is also a~nonlinear transformation, which relates the
solution $w(x,t)$ to the solution $\phi(x,t)$ of the NDE
\eqref{FBE:03B}
\begin{gather}
\label{FBE:62}
\phi(x,t)=-2\alpha\log{w(x,t)}.
\end{gather}

We see that the above two transformations are the cases of one and
the same transformation which is as follows
\begin{gather}
\label{FBE:63}
\phi(x,t)=-2\alpha\,{}_{a}D^{p}_x\log{w(x,t)},
\end{gather}
for $p=0$~\eqref{FBE:62} and~$p=1$~\eqref{FBE:61}, and~the
transformation~\eqref{FBE:63} itself relates the dif\/fusion
equation solution $w(x,t)$ to the FBENN~\eqref{fbe:0} solution
$\phi(x,t)$.

It is more expedient, however, particularly in applications, to
use the fractional generalization of the Hopf--Cole
transformation
\begin{gather}
\label{FBE:64}
\phi(x,t)=-2\alpha\,{}_{a}D^{p}_x\log{[b+w(x,t)]},
\qquad b\in\mathbb{R}.
\end{gather}

Indeed, by substituting expression~\eqref{FBE:64} into the initial
FBENN~\eqref{fbe:0}, we obtain the equation{\samepage
\begin{gather*}
\,_{a}D^{p}_x\left[\frac{w_t-\alpha w_{xx}}{b+w}\right]=0,
\end{gather*}
which turns into an identity when $w(x,t)$ is the solution of the DE.}

The point is that by changing the scale and~the variables we can
always obtain the solution of the dif\/fusion equation for
$||\,w(x,t)||\ll1$.
In this case, the generalization of the fractal
Hopf--Cole transformation~\eqref{FBE:64} becomes even simpler
\begin{gather}
\label{FBE:65}
\phi(x,t)=-2\alpha\,{}_{a}D^{p}_xw(x,t).
\end{gather}

From the formula~\eqref{FBE:65} in particular it follows that if
the dif\/fusion equation solution asymptotically approaches zero,
for instance, for $w(x,t)\mathop{\rightarrow}x^{-q}$ for $x
\to+\infty$ $(q>0)$, then the FBENN solution also
approaches zero: $w(x,t)\mathop{\rightarrow}
x^{-(p+q)}$ for $x\to+\infty$, $(q>0)$.
In the case of
the exponential asymptotic, more convenient is the expression~\eqref{FBE:63}.

Let $w(x,t):x\in[0,+\infty)\cup t\in
[0,+\infty)$.
The fractional derivative is the right Caputo
derivative in the Weyl sense.
In the case of the solution of the
dif\/fusion equation
\begin{gather*}
w(x,t)=\exp{\left(-\frac{cx}{2\alpha}+\frac{c^2t}{4\alpha}-b\right)},
\end{gather*}
the solution of the FBENN is
\begin{gather}
\label{27t}
\phi(x,t)=-2\alpha
\log{\left[a+\exp{\left(-\frac{cx}{2\alpha}+\frac{c^2t}{4\alpha}-b\right)}\right]},\qquad p=0,
\\
\label{28t}
\phi(x,t)=
\frac{c}{a+\exp{\left(-\frac{cx}{2\alpha}+\frac{c^2t}{4\alpha}-b\right)}},\qquad p=1,
\\
\label{29t}
\phi(x,t)=-2\alpha\,{}_{+\infty}D^{p}_{x-}
\log{\left[a+\exp{\left(-\frac{cx}{2\alpha}+\frac{c^2t}{4\alpha}-b\right)}\right]},\qquad0<p<1.
\end{gather}

Note that solution~\eqref{29t} continuously transforms from
solution~\eqref{27t} into~\eqref{28t} when the parameter $p$ runs
from $p=0$ to $p=1$.

It is the fractional generalization of the Hopf--Cole
transformation~\eqref{FBE:61} that has been used to derive the
solution~\eqref{29t} which interrelates the dif\/fusion and~the
FBENN solutions.

Thus, if $T^{(p)}$ is a~fractional generalization of the
Hopf--Cole transformation, then the interrelation among the
nonlinear dif\/fusion equation (NDE), the Burgers equation (BE), the
simple dif\/fusion equation (DE) and~the nonlinear nonlocal
dif\/fusion equation can be graphically shown as follows
$$
\begin{gathered}
\hspace*{27mm}\textrm{DE}\phantom{AAAAAAl}\nonumber \\
\hspace*{14mm}\swarrow {}^{T^{(0)}} \ \ \ \  \ \downarrow {}^{T^{(p)}} \ \ \ \ \searrow
\!\!\!\!{}^{T^{(1)}}\phantom{A}
\nonumber \\
\phantom{AA}\textrm{NDE} \xleftarrow [p\rightarrow 0]{}
\textrm{FBENN} \xrightarrow[p\rightarrow 1]{} \textrm{BE}
\nonumber
\end{gathered}
$$

\section{The interrelation of the solutions}
\label{section6}

The existence of the fractional generalization of the Hopf--Cole
transformation may produce the impression that the properties of
the FBENN solutions can be reduced to the nonlocally transform
solution of the dif\/fusion equation.
Below, we shall show this is
not the case.

It is well-known that any solution $w(x,t)$ of the DE
$w_t=\alpha w_{xx}$ under the simple Hopf--Cole
transformation $\phi(x,t)=-2\alpha w_x(x,t)/w(x,t)$ turns into
the solution $\phi(x,t)$ of the BE $\phi_t+\phi\phi_x-\alpha
\phi_{xx}=0$.
The reverse is not true, as the solution
$\phi(x,t)$ does not obey the DE: if we substitute $\phi(x,t)$
into the BE, we obtain the equation which will be more general
than the~DE
\begin{gather*}
w_t-\alpha w_{xx}=f(t)w,
\end{gather*}
here $f(t)$ is any time function.

Thus, there exist the BE solutions that can be expressed through
solutions of the simple DE, and~the proper solutions that are
devoid of such representation.

In a~similar way, we can classify also the NDE solutions.
Some of them are related to DE solutions through the transformation
$\phi(x,t)=-2\alpha\log{w(x,t)}$.
However, there exist also the
proper solutions that are not related to DE solutions.

It is important that the interrelation of the solutions in the
case of the travelling wave equation~\eqref{FBE:14} can be
generalized and~is valid for the NDE, FBENN and~BE cases
\begin{gather*}
\phi(x,t)={}_{a}D^{p}_x\phi^{(0)}(x,t),\qquad
\phi(x,t)={}_{a}D^{p-1}_x\phi^{(1)}(x,t),
\end{gather*}
here, as above, $\phi^{(0)}(x,t)$ ($\phi^{(1)}(x,t)$) is the
NDE (BE) solution for $p=0$ ($p=1$).
This means that the
solutions of all these three equations, irrespectively of the
fractional Hopf--Cole transformation, are interrelated.
Hence, in
particular, it follows that
\begin{gather*}
\phi^{(1)}(x,t)=\partial_x\phi^{(0)}(x,t).
\end{gather*}

This interrelation is, in a~sense, more general as is valid
for both the solutions related to the DE solutions and~proper
solutions.
If we substitute the proper solution $\phi^{(1)}(x,t)$
of the BE, we will be able to restore the proper solutions of the
NDE.

Let us consider the initial evolution equation~\eqref{fbe:0}
together with two limit cases at $p\rightarrow0$ and
$p\rightarrow1$
\begin{gather}
\label{FBE:03}
\phi_{t}+\frac{1}{2} \phi^{2}_{x}-\alpha
\phi_{xx}=0,
\\
\label{FBE:04}
\phi_{t}+\phi\phi_{x}-\alpha\phi_{xx}=0.
\end{gather}

For $p=0$ we have a~dif\/fusion equation with quadratic
nonlinearity or just the nonlinear dif\/fusion equation (NDE)
\eqref{FBE:03}, for $p=1$~-- the Burgers equation (BE)
\eqref{FBE:04}, for $0<p<1$~-- the integro-dif\/ferential equation
FBENN~\eqref{fbe:0} which we may regard as an ``intermediate''
evolution equation whose solutions turn into the solution of
equation~\eqref{FBE:03} or~\eqref{FBE:04} depending on $\mathop
{\lim}\limits_{p\to\,0}\phi^{(p)} ({x, t} )$ or
$\mathop{\lim}\limits_{p\to1}\phi^{(p)} ({x, t}
 )$.

The solutions to the Cauchy problem~\eqref{fbe:0} have to be
understood in some weak sense; there are several options
reported, e.g., in~\cite{BilerFunakiWoyczynski1998,BilerKarchWoyczynski2001}.
In the context of the present study, let us just say that under the weak
solution~\eqref{fbe:0} we mean the solution of the integral equation
\begin{gather}
\label{FBE:05}
\phi ({x, t} )=e^{\alpha t\partial_{xx}}\phi_0
-\frac{1}{2}\int_0^t{e^{\alpha ({t-\tau} )
\partial_{xx}}}\cdot {}_{a}D^p_x\big({{}_{a}D^{1-p}_x\phi}\big)^2 ({x, t} )d\tau
\end{gather}
motivated by the classical Duhamel formula.

It is clear now that in the same manner we can also classify the
nonlocal solutions of the FBENN.
Some of them, through the
fractional Hopf--Cole transformation, are related to the simple
DE.
However, there are also a~number of proper solutions which
have no such relation to the BE.

\section{Evolution of the initial conditions}
\label{section7}

Let the initial conditions for FBENN~\eqref{fbe:0} and~for the
dif\/fusion equation be related by the expression
\begin{gather*}
w_0(x)=e^{-\frac{1}{2\alpha}\,{}_{a}I^{p}_x\phi_0(x)}.
\end{gather*}
This allows us to express the solution $\phi(x,t)$ of the FBENN
through the initial condition $\phi_0(x)$ and~the general form
of the solution of the dif\/fusion equation:
\begin{gather}
\label{FBE:42}
\phi(x,t)=-2\alpha\,_{a}D^p_x\left[\log{\left(1+
\frac{1}{\sqrt{4\pi\alpha t}}
\int_{-\infty}^{+\infty}e^{-\frac{|x-y|^2}{4\alpha t}
-\frac{1}{2\alpha}\,{}_{a}I^{p}_y\phi_0(y)}dy
\right)}\right].
\end{gather}

The expression of the solution of the FBENN~\eqref{fbe:0} in the general
form~\eqref{FBE:42} allows us to analyze the time evolution of the
nonlocal solution.
Indeed, let $w(x,t)$ be the solution of the
dif\/fusion equation $w_t=\alpha w_{xx}$.
Then the solution of
\eqref{FBE:42} can be expressed in the form $\phi(x,t)=-2\alpha
\,{}_{a}D^{p}_x\log{(1+w(x,t))}$.

As follows from the FBENN~\eqref{fbe:0}, depending on the values
of the parameter $p$ we deal with not one but with an inf\/inite
number or an hierarchy of integro-dif\/ferential equations.
One of
the most important properties of the FBENN~\eqref{fbe:0} is
interrelation between nonlinearity and~nonlocality: for the
fractional value of the parameter~$p$ we have the
nonlinear-nonlocal and~for the integer positive~$p$ only a
nonlinear generalization of the Burgers equation.
In this
hierarchy, due to the substitution $\phi(x,t)\rightarrow
{}_{a}D^{q-p}_x\phi(x,t)$, the low-order equations
turn into the higher-order ones, but in the inverse direction this
transformation is multivalued.

\section{Interaction of nonlinear and~nonlocal perturbations}
\label{section8}

The relation~\eqref{FBE:42} between the solutions of the FBENN and
dif\/fusion equation allows to consider an interaction of nonlocal
and nonlinear perturbations.
Two or more perturbations moving with
a dif\/ferent velocity can overtake each other or f\/low together into
a new intensive perturbation.
The FBENN also describes the
interaction process of two or more moving nonlocal perturbations.
The principle of superposition is not valid for the nonlinear
FBENN, but it is valid for the linear dif\/fusion equation.
The
fractional Hopf--Cole transformation~\eqref{FBE:42} interrelates
the solutions of the nonlocal and~nonlinear FBENN and~of the
linear dif\/fusion equation.
Thus, if $w_i(x,t)$ are the
solutions of the dif\/fusion equation, then $\phi(x,t)=-2\alpha\,
{}_{a}D^p_x(\log{\sum{w_i}})$ are the solutions of the
FBENN.

For instance, for two solutions of the dif\/fusion equation in the
form
\begin{gather*}
w_i(x,t)=a_i\exp{\left(-\frac{c_ix}{2\alpha}+\frac{c^2_it}{4\alpha}-b_i\right)},\qquad i=1,2,
\end{gather*}
we obtain a~nonlocal and~nonlinear interaction of these
perturbations:
\begin{gather*}
-\frac{\phi(x,t)}{2\alpha}=
\begin{cases}
\log{\big(w_1+w_2\big)},&p=0,
\vspace{1mm}\\
\displaystyle
\frac{1}{\Gamma(1-p)}\frac{d}{dx}\int^x_a\frac{\log{[w_1(\xi,t)+w_2(\xi,t)]}}{(x-\xi)^{p}}\,d\xi,\qquad&0<p<1,
\vspace{1mm}\\
\displaystyle
\frac{c_1w_1+c_2w_2}{w_1+w_2},& p=1.
\end{cases}
\end{gather*}

\section{The conservation laws}
\label{section9}

In the case of the BE~\eqref{FBE:04}, for $x\in E^1$, $\forall\,
t>0$, $\phi(\pm\infty,t)=\phi_x(\pm\infty,t)=0$, we have a
conservation value, or the time invariant value $\operatorname{inv}$, of
\begin{gather*}
I^{(1)}=\int_{-\infty}^{+\infty}\phi(x,t)\,dx=\operatorname{inv},
\end{gather*}
since
\begin{gather*}
\frac{\partial I^{(1)}}{\partial t}
=\int_{-\infty}^{+\infty}\left[\alpha\phi_x-\frac{1}{2}\phi^2\right]_xdx
=\left(\alpha\phi_x-\frac{1}{2}\phi^2\right)\bigg|^{+\infty}_{-\infty}=0.
\end{gather*}

In the applications, this conservation law is called the ``mass'' 
conservation law, because if $\phi(x,t)$ can be a~one-dimensional
density or a~gradient of any physical, chemical or biological
magnitude.
Then, the $I^{(1)}$ corresponds to its conservation.

In the case of the NDE~\eqref{FBE:13}, if $\forall\,t>0$,
$\phi_x(\pm\infty,t)=\phi_{xx}(\pm\infty,t)=0$, we again
deal with the conservation value
\begin{gather}
\label{FBE:83}
I^{(0)}=\phi(+\infty,t)-\phi(-\infty,t)=\operatorname{inv},
\end{gather}
since by applying the derivative $\partial_x$ to the
evolutionary equation~\eqref{FBE:13} followed by integration we obtain
\begin{gather*}
\frac{\partial I^{(0)}}{\partial t}=
\frac{\partial}{\partial t}\int_{-\infty}^{+\infty}\phi_x\,dx=
\left(\alpha\phi_{xx}-\frac{1}{2}\phi_x^2\right)\bigg|^{+\infty}_{-\infty}=0.
\end{gather*}

This conservation value shows that the dif\/ference in asymptotic
values for any time moment remains unchanged.
If, for instance, we
deal with the evolution of the potentials, the conservation value
$I^{(0)}$~\eqref{FBE:83} shows that the dif\/ference of potentials
for $x\rightarrow\pm\infty$ does not change.
In the case of the FBENN~\eqref{FBE:14}, we again deal with a
conservation value if $\forall\,t>0$,
$_{a}D^{2-p}_x\phi(\pm\infty,t)={}_{a}D^{1-p}_x\phi(\pm\infty,t)=0$
\begin{gather}
\label{FBE:85}
I^{(p)}=\int_{-\infty}^{+\infty}{}_{a}D^{1-p}_x\phi(x,t)\,dx=\operatorname{inv},
\end{gather}
because
\begin{gather*}
\frac{\partial I^{(p)}}{\partial t}=
\int_{-\infty}^{+\infty}\left[\alpha\,{}_{a}D^{2-p}_x\phi-\frac{1}{2}\big({}_{a}D^{1-p}_x\phi\big)^2\right]_xdx=
\left[\alpha\,{}_{a}D^{2-p}_x\phi-\frac{1}{2}\big({}_{a}D^{1-p}_x\phi\big)^2\right]\bigg|^{+\infty}_{-\infty}=0.
\end{gather*}

Even this simple example highlights two important properties of
the nonlocal conservation law~\eqref{FBE:85}: it interrelates the
conservation values of two dif\/ferent dynamical systems, which can
be of dif\/ferent mathematical nature (e.g., in our case
these values are integral and~discrete).

Note that in the ``common'' case of the nonlocal BE
\begin{gather*}
\phi_t+\phi\phi_x-\alpha\,{}_{a}D^{2-p}_x\phi=0,
\end{gather*}
an analogous conservation integral exists at other asymptotic values
\begin{gather*}
\phi(\pm\infty,t)={}_{a}D^{1-p}_x\phi(\pm\infty,t)=0.
\end{gather*}

At this point, it is not dif\/f\/icult to characterize the ``mass''
conservation law of the nonlinear nonlocal evolution equation
\begin{gather}
\label{FBE:88}
\phi_t+\frac{1}{2}\,{}_{a}D^{p}_x\big({}_{a}D^{1-p}_x\phi\big)^2-\alpha\,_{a}D^{2-q}_x\phi=0.
\end{gather}
The magnitude $I^{(p,q)}$ is the invariant of the evolution
equation~\eqref{FBE:88}:
\begin{gather*}
I^{(p,q)}=\int_{-\infty}^{+\infty}{} _{a}D^{1-p}_x\phi(x,t)\,dx=\operatorname{inv}
\end{gather*}
for $_{a}D^{1-p}_x\phi(\pm
\infty,t)={}_{a}D^{2-(p+q)}_x\phi(\pm\infty,t)=0$.
However,
in this case the integrability is sacrif\/iced, and~the fractional
generalization of the Hopf--Cole transformation does not exist.

The ``energy'' of travelling excitation
\begin{gather}
\label{FBE:810}
K=\frac{1}{2}\int_{-\infty}^{+\infty}\big({}_{a}D^{1-p}_x\phi\big)^2\,dx,
\end{gather}
if $\forall\, t>0$, ${}_{a}D^{2-p}_x\phi(\pm\infty,t)={}_{a}D^{1-p}_x\phi(\pm\infty,t)=0$,
as in the case of the BE, is not unchangeable, but it is constantly decreasing:
\begin{gather*}
\frac{d}{dt} \frac{1}{2}\int_{-\infty}^{+\infty}\big({}_{a}D^{1-p}_x\phi\big)^2\,dx=
-\frac{1}{3}\big({}_{a}D^{1-p}_x\phi\big)^3\bigg|^{+\infty}_{-\infty}
+\alpha\big({}_{a}D^{1-p}_x\phi\big)\big({}_{a}D^{2-p}_x\phi\big)\bigg|^{+\infty}_{-\infty}
\\
\qquad
{}-\alpha\int_{-\infty}^{+\infty}\big({}_{a}D^{2-p}_x\phi\big)\cdot
\big({}_{a}D^{2-p}_x\phi\big)\,dx
=-\alpha\int_{-\infty}^{+\infty}\big({}_{a}D^{2-p}_x\phi\big)^2\,dx<0.
\end{gather*}

Like in the case of the nonlocal ``mass'' conservation law, the
``energy'' $K$~\eqref{FBE:810} links the energy of the travelling
excitation in the case of the BE $(p=1)$ and the one-dimensional
density of energy in the case of the NDE $(p=0)$.

\section{Symmetries}
\label{section10}

The conservation laws are obviously related to the group of the FBENN
automorphisms.
In the case of the usual BE, there is a~large
Lee symmetry group of point transformations, which contains the
Galilei, dilaton, and~projective transformations and~is generated
by operators
\begin{gather*}
D=2t\partial_t+x\partial_x-\phi\partial_{\phi},
\qquad
K=t^2\partial_t+tx\partial_x-\left(tu+\tfrac{1}{2}x\right)
\partial_{\phi},
\\
P_1=\partial_x,\qquad P_2=\partial_t,\qquad
B=2t\partial_x-\partial_{\phi}.
\end{gather*}
Symmetries in the discrete BE have been studied in~\cite{HLW1999}.

An important and~in the general case uninvestigated problem is to f\/ind
the nonlocal and~nonclassical symmetries of the FBENN.
Some general
aspects were already presented in~\cite{AbrahamGuo1993}.
Promi\-sing seems an attempt of computer symmetry analysis, as was
done for the nonlinear heat equation in~\cite{Clarcson_Mansfield1993}.

Note here one property that allows us to get new solutions of the
FBENN.
Let $v(x,t)$ be a~known solution of the FBENN
\eqref{fbe:0}, and~$u(x,t)$ is the solution of the linear equation
\begin{gather}
\label{FBE:91}
u_t+\big({}_{a}D^{1-p}_x v\big) u_x-\alpha u_{xx}=0.
\end{gather}
Then
\begin{gather*}
\phi(x,t)=-2\alpha\,{}_{a}D^{p}_x\log{u}+v
\end{gather*}
is a~new solution of the FBENN~\eqref{fbe:0}.
At an integer value
of $p$ we have the local equation~\eqref{FBE:91}; in particular,
for $p=1$ we obtain a~new solution of the classical Burgers
equation.

Note here that the existence of transformation $T$ expressed by
relation~\eqref{FBE:64} allows to solve the problem
of the FBENN~\eqref{fbe:0} symmetry group.
If $G_1$ is a
symmetry group of the dif\/fusion equation, then $G=TG_1T^{-1}$
is a~symmetry group of the FBENN~\eqref{fbe:0}.

Here, we shall note an important feature of the FBENN, related to the FBENN symmetry.
Let us apply the operator $\partial_x$ to the FBENN and~take the operator $_{a}D^{p}_x$ beyond the
brackets
\begin{gather*}
{}_{a}D^{p}_x\Big[{}_{a}D^{1-p}_x\phi_t+
\big({}_{a}D^{1-p}_x\phi\big)\big({}_{a}D^{1-p}_x\phi\big)_x-
\alpha
\big({}_{a}D^{1-p}_x\phi\big)_{xx}\Big]=0.
\end{gather*}
On substituting the variables while
\begin{gather*}
x\mapsto y\pm A(t+B),\qquad{}_{a}D^{p}_x\phi(x,t)
\mapsto {}_{a}D^{p}_y\phi(y,t)\pm A,
\end{gather*}
where $A$ and~$B$ are constants, the evolution FBENN will not
change its form.
This way of generating the new solutions is
particularly ef\/f\/icient not in the case of travelling excitations
when $\phi=\phi(x-ut)$, but in the general case when
$\phi=\phi(x,t)$.

One of the ways of constructing exact solutions of some of the
nonlinear equations consists in f\/inding the corresponding
B\"{a}klund transformations.
The B\"{a}klund transformations have
been found for the majority of the nonlinear equations that are
integrated by the inverse scattering method~\cite{Zakharov1984,ZakharovFokas1993}.

It is possible to show that the generalization of the Hopf--Cole
transformation $T^{(p)}$ is a~separate case of the B\"{a}klund
transformation
\begin{gather}\label{FBE:96}
w_x+\frac{1}{2\alpha}(b+w)\,{}_{a}D^{1-p}_x\phi=0,
\qquad
w_t+\frac{1}{2}\big[(b+w)\,{}_{a}D^{1-p}_x\phi\big]_x=0.
\end{gather}

On removing $\phi(x,t)$ from the system we obtain the dif\/fusion
equation $w_t=\alpha w_{xx}$.

To obtain the FBENN from the system~\eqref{FBE:96}, we shall
perform the following procedure.
From the dif\/fusion equation for
the function $w(x,t)$ and~the identity
$[w_t/(b+w)]_x=[w_x/(b+w)]_t$ follows the condition
\begin{gather}
\label{FBE:97}
\left(\frac{w_x}{b+w}\right)_t=\alpha\left(\frac{w_{xx}}{b+w}\right)_x,
\end{gather}
and from the f\/irst equation of the system~\eqref{FBE:96} it
follows that
\begin{gather}
\label{FBE:98}
\frac{w_x}{b+w}=-\frac{1}{2\alpha}\,{}_{a}D^{1-p}_x\phi,
\qquad\mbox{and}\qquad
\frac{w_{xx}}{b+w}=-\frac{\big({}_{a}D^{1-p}_x\phi\big)_x}{2\alpha}+
\frac{\big({}_{a}D^{1-p}_x\phi\big)^{2}}{4\alpha^2}.
\end{gather}

Substitution of the expressions~\eqref{FBE:98} in the condition
\eqref{FBE:97} gives us the FBENN~\eqref{fbe:0}
\begin{gather*}
{}_{a}D^{1-p}_x\left[
\phi_t+\frac{1}{2}\,{}_{a}D^p_x\big({}_{a}D^{1-p}_x\phi\big)^2-\alpha\phi_{xx}
\right]=0.
\end{gather*}

The f\/irst equation of the system~\eqref{FBE:96} is a~fractional
generalization of the Hopf--Cole transformation
\begin{gather}
\label{FBE:99}
\phi(x,t)=-2\alpha\,{}_{a}D^{p}_x\log{[b+w(x,t)]}.
\end{gather}

On substituting expression~\eqref{FBE:99} into the FBENN we
obtain the following equation
\begin{gather}
\label{FBE:910}
{}_{a}D^{p}_x\left[\frac{1}{b+w} (w_t-\alpha w_{xx} )\right]=0,
\end{gather}
from which it follows that any solution of the dif\/fusion equation
can be transformed, through the generalization~\eqref{FBE:99} of
the fractional Hopf--Cole transformation, into an FBENN solution.
The reverse implication would be incorrect, as from equation~\eqref{FBE:910} follows a~more general dif\/fusion equation.
For instance,
\begin{gather*}
w_t-\alpha w_{xx}= (C_1+f(t) )(b+w).
\end{gather*}

It is because of this above-mentioned ambiguity of the fractional
generalization of the Hopf--Cole transformation that the FBENN
solutions are dif\/ferentiated into those related to the dif\/fusion
equation and~proper ones.

\section{The asymptotic form of solutions}
\label{section11}

A very convenient dimensionless quantity which is used in the
nonlinear BE is the Reynolds number.
This number is just a~ratio
of the nonlinear and~the dissipative terms: ${\textsf{Re}}\sim
\phi\phi_{x}/\alpha\phi_{xx}$.
In the case when
${\textsf{Re}}\ll1$, the inf\/luence of the nonlinear term is
negligible, but for ${\textsf{Re}}\gg1$ this term plays a~crucial
role and~leads to the nonlinear Riemann equation which describes the
simplest type of the shock waves\footnote{This equation has some
names: Hopf, weak shock waves, Riemann--Hopf.
We are following
reference~\cite{Riemann1860}.}.

For the FBENN~\eqref{fbe:0} we may introduce a~dimensionless
generalization of the Reynolds number, which at the $p=1$ coincides
with the classical Reynolds number for the BE
\begin{gather*}
{\textsf{Re}}\sim\frac{
{}_{a}D^p_x\big({}_{a}D^{1-p}_x\phi\big)^2}
{\alpha\phi_{xx}}\sim\frac{\phi x^p}{\alpha\lambda^{1+p}},
\end{gather*}
where $\lambda$ is the characteristic parameter of length in the
model.

Depending on the value of this number, we obtain two limit cases
of the FBENN
\begin{gather}\label{FBE:102}
\phi_t-\alpha\phi_{xx}=0,\qquad{\textsf{Re}}\ll1,
\\
\phi_t+\frac{1}{2}\,{}_{a}D^p_x\big({}_{a}D^{1-p}_x\phi
\big)^2=0,\qquad{\textsf{Re}}\gg1.\label{FBE:102a}
\end{gather}

The equation~\eqref{FBE:102} is just a~dif\/fusion
equation, and~the  equation~\eqref{FBE:102a} we will call the fractional Riemann
equation.
In the case of ${\textsf{Re}}\gg1$, depending on the
value of the order of the fractional derivative $p$, we have three
evolution equations
\begin{gather*}
\phi_t+\frac{1}{2} \phi^{2}_{x}=0,\qquad p=0,
\\
\phi_t+\frac{1}{2}\,{}_{a}D^p_x\big({}_{a}D^{1-p}_x\phi\big)^2=0,\qquad 0<p<1,
\\
\phi_t+\phi\phi_{x}=0,\qquad p=1.
\end{gather*}

The ``mass'' conservation law predetermines the asymptotic form of
the FBENN solution.
To conf\/irm such a~result, let us consider
some estimates.
From the general form of the FBENN solution
\eqref{FBE:42} it follows that the limit $t\rightarrow+\infty$
corresponds to a~rather low value of the parameter~$\alpha$.
At a
low~$\alpha$, to calculate the values of the corresponding
integrals we can apply the saddle point approximation.

The critical point $y_0$ can be determined from the equation
$\frac{y_0-x}{t}+{}_{a}D^{1-p}_{y}\phi_0(y_0)=0$.
Then the asymptotic expression of the FBENN solution acquires the form
\begin{gather}
\label{FBE:105}
\phi(x,t)={}_{a}D^p_x\left[\frac{(x-y_0)^2}{2t}+C\right]
\sim
\frac{(x-a)^{2-p}}{t\Gamma(3-p)}+C_1(x-a)^{-p}.
\end{gather}
For $x\rightarrow+\infty$ and~$0<p<1$, the solution
$\phi(x,t)\rightarrow(x-a)^{2-p}/[t\Gamma(3-p)]$.
Thus, we obtain
a~po\-wer-deformed perturbation of the usual solution of the Burgers equation.
Note here that these estimates are valid not only in the
environment of the meaning $p=1$, but also for any $p\in
\mathbb{R}$.

We have to show the region of the validity of solution~\eqref{FBE:105}.
In the limit case, the integral in the expression
of the mass conservation law diverges.
Therefore, for $x>x_0$
the solution $\phi(x)\equiv0$.
To determine the value $x_0$,
we insert the asymptotic form of solution~\eqref{FBE:105} in the
expression of the mass conservation law~\eqref{FBE:85}.
This means
that $x^2_0/2t\sim I$.
Thus, the maximum meaning of the solution
\begin{gather*}
\phi_{\max}(x,t)\sim\frac{I^{1-\frac{p}{2}}}
{(2t)^{\frac{p}{2}}}\qquad\mbox{and}\qquad x_0\sim\sqrt{2It}.
\end{gather*}

\section{The supersymmetric nonlinear nonlocal    dif\/fusion\\ evolution equation}
\label{section12}

We shall show that the FBENN has a~supersymmtric generalization.
Let the superf\/ield $\chi=\theta\,{}_{a}D^{1-p}_x\phi+\psi$
unite two f\/ields of dif\/ferent properties: the ``bosonic'' f\/ield
$\phi(x,t)$ and~its spinor superpartner $\psi(x,t)$; $\theta$ is
the constant Majorana spinor.
The transformations of the f\/ields
$\phi,\psi$ are nonlocal because of the fractional derivatives
${}_{a}D^{p}_xf(x)$
\begin{gather}\label{FBE:111}
\delta_{\eta}\psi=\eta\,{}_{a}D^{1-p}_x\phi,
\qquad
\delta_{\eta}\,{}_aD^{1-p}_x\phi=\eta\psi_x.
\end{gather}
However, the commutator of the two transformations~\eqref{FBE:111}
is a~spatial translation
\begin{gather*}
\left[\delta_{\eta},\delta_{\xi}\right]=2\xi\eta\partial_x.
\end{gather*}

The supersymmtric equation
\begin{gather}
\label{FBE:113}
\chi_t=\left(\chi_x+\frac{1}{2}\chi\mathcal{D}\chi\right)_x
\end{gather}
(here $\mathcal{D}=\theta\partial_x+\partial_{\theta}$ is a
supersymmetric derivative) is a~system of two evolutionary
equations
\begin{gather}
\label{FBE:114}
\psi_t=\psi_{xx}+\frac{1}{2}\big({}_{a}D^{1-p}_x\phi\cdot\psi\big)_x,
\qquad
\phi_t=\phi_{xx}+\frac{1}{2}\,{}_{a}D^{p}_x\left[\big({}_{a}
D^{1-p}_x\phi\big)^2-\psi\psi_x\right]
\end{gather}
which are invariant in respect of supertransformations
\eqref{FBE:111}.
In the general case, the system~\eqref{FBE:114}
is a~system of two nonlinear nonlocal evolution equations, which
becomes local when
\begin{gather*}
p=0\qquad\begin{cases}
\psi_t=\psi_{xx}+\frac{1}{2} (\phi_x\psi )_x,
\\
\phi_t=\phi_{xx}+\frac{1}{2} \big(\phi^2_x-\psi\psi_x \big),
\end{cases}
\\
p=1\qquad\begin{cases}
\psi_t=\psi_{xx}+\frac{1}{2} (\phi\psi )_x,
\\
\phi_t=\phi_{xx}+\frac{1}{2}\big(\phi^2-\psi\psi_x\big)_x.
\end{cases}
\end{gather*}

The supersymmetric equation~\eqref{FBE:113} and~the corresponding
system of equations~\eqref{FBE:114} unite two f\/ields of dif\/ferent
nature, and~only one of them is nonlocal.

Here, the main point is a~general note related to the application
of the nonlocal systems.
Suppose the case when a~dynamic system is
characterized by two interacting f\/ields, of them one, for
instance, the ``fermionic'' f\/ield $\psi(x,t)$, is measured in
the course of experiment, whereas the other, the ``bosonic'' f\/ield
$\phi(x,t)$, is assessed only phenomenologically.
Actually, such assessment in the class of \emph{local} evolution equations may result in a
qualitatively erroneous mathematical model of a~dynamic system.

\section{Conclusions}
\label{section13}

It is important to note that the inf\/luence of nonlocality can be
arbitrarily great.
Therefore, we do not describe nonlocality by an
additional term in the Burgers equation.

Let us remind here that the classical Burgers equation belongs to a
unique group of the three completely integrable second-order PDEs.
I suggest that the FBENN also belongs to a~unique group of the
completely integrable nonlocal PDEs of the fractional order.

The fractional dif\/fusion process is related to the non-Gauss
statistics; this results in slow dif\/fusion correlators $\langle
(\Delta x)^2\rangle\propto Dt^{\gamma}$ with $\gamma\neq1$, and
$D$ is a~generalized dif\/fusion coef\/f\/icient of the dimension
$L^2/T^{\gamma}$.
In our case, the FBENN is related to the
so-called L\'{e}vi statistics~\cite{ShlesingerZaslavskyFrisch}; at
the same time, the initial Burgers equation as well as the
dif\/fusion equation are related to the usual Gauss statistics.

From our point of view, the FBENN has at least two important advantages:
\begin{itemize}\itemsep=0pt
\item[$i)$] the inf\/luence of nonlocality is not assumed to be
insignif\/icant;
\item[$ii)$] the relation of the FBENN to the usual dif\/fusion equation allows a~lot of analytical solutions of the FBENN.
\end{itemize}

Besides, despite the nonlocality in the proposed nonlinear and nonlocal FBENN,
\begin{itemize}\itemsep=0pt
\item  space-localized solutions are possible;
\item  nonlocal perturbations in a~system described by the FBENN can interact;
\item  the Reynolds number is a~universal dimensionless parameter for both the local and~nonlocal Burgers equations;
\item  there are nonlocal analogies of both the momentum conservation law and kinetic energy dissipation.
\end{itemize}

In some f\/ields of physics, we need the vectorial form of the
Burgers equation, e.g., in astrophysics to describe the large-scale
structure of the Universe~\cite{Zeldovich1970,Miskinis2000,Miskinis2003}.
In such cases, the vectorial FBENN can be proposed
\begin{gather*}
{\bf\phi}_t+\frac{1}{2}\,_{a}{D}^p_x\big({}_{a}{D}^{1-p}_x\phi\big)^2
-\alpha{\bf\nabla}({\bf\nabla\phi})=0,
\end{gather*}
where ${\bf\phi}=(\phi^{1},\ldots,\phi^{n})\in {\mathbb R}^n$, ${}_{a}{D}^p_x$
is the fractional generalization of the gradient operator~${\bf\nabla}$.

Note also that for $\alpha=0$ from FBENN~\eqref{fbe:1} follows the
fractional generalization of the Riemann equation, which also has
numerous applications.

The proposed FBENN, because of its general character, allows a
wide range of applications.
Actually, we may try to introduce the
nonlocal generalization in almost all f\/ields where the BE is
applied.
These are the nonlocal ef\/fects in shock wave propagation
in acoustics, the ef\/fective model of the process of nonlinear heat
distribution in the environment in the presence of heat sources
and sinks, the Kardar--Parisi--Zhang (KPZ) e\-qua\-ti\-on in the
crystal growth phenomena in (1+1)-di\-men\-si\-ons~\cite{KardarParisiZhang1986},
the nonlinear dynamics of moving lines \cite{HwaKardar1992},
formation of galaxies~\cite{ShandarinZeldovich1989,Miskinis2000}, the behavior of the magnetic f\/lux line
in su\-per\-con\-duc\-tors~\cite{Hwa1992}, and~spin glasses
\cite{FisherHuse1991}, as well as numerous examples of the
application of the usual Burgers equation, presented in the above-mentioned monographs
\cite{GurbatovMalakhovSaichev1991,Woyczynski1998}.

\appendix
\section{The normalized form of equations}
\label{appendixA}

In applications, usually considered are both the non-normalized form of
the dif\/fusion equation with quadratic nonlinearity
\begin{gather}
\label{FBE:A1}
\phi_{t}+\alpha\phi^{2}_{x}-\beta\phi_{xx}=0
\end{gather}
and the Burgers equation
\begin{gather}
\label{FBE:A2}
\phi_{t}+\alpha\phi\phi_{x}-\beta\phi_{xx}=0,
\end{gather}
in which the concrete sense of the coef\/f\/icients~$\alpha$ and~$\beta$
depends on the content of a~model under description.
In the
dimensional form, $\beta$ is normally related to the coef\/f\/icient
of dif\/fusion, whereas in the dimensionless form it relates to the
inverse Reynolds number $\textsf{Re}^{-1}$.

Equations~\eqref{FBE:A1} and~\eqref{FBE:A2} in the non-normalized
form, as well as the non-normalized equation FBENN~\eqref{fbe:0}
\begin{gather*}
\phi_t+\alpha\,{}_{a}D^p_x\big({}_{a}D^{1-p}_x\phi\big)^2-\beta\phi_{xx}=0,
\end{gather*}
on performing the transformations
\begin{gather*}
\phi(x,t)\mapsto\frac{\beta}{\alpha} \phi(x,t),\qquad t
\mapsto\frac{1}{\beta} t,
\end{gather*}
are reduced to the normalized form
\begin{gather*}
\phi_{t}+\phi^{2}_{x}+\phi_{xx}=0,
\qquad
\phi_t+{}_{a}D^p_x\big({}_{a}D^{1-p}_x\phi\big)^2+\phi_{xx}=0,
\qquad
\phi_{t}+\phi\phi_{x}+\phi_{xx}=0.
\end{gather*}

Note that in the case of the Burgers equation (for $p=1$) there exists
a special transformation related to the changed scale of the
independent variables without changing the function,
\begin{gather*}
x\mapsto\frac{\beta}{\alpha}\,x,\qquad t\mapsto\frac{\beta}{\alpha^2}\,t,
\end{gather*}
which allows also a~reduction of the Burgers equation to a
dimensionless form~\eqref{FBE:04A}.
Thus, now we may apply any form
of the FBENN or its limit cases the BE or the NDE depending on our requirements.

\section{The fractional calculus}
\label{appendixB}

The left-side Riemann--Liouville fractional derivative of the order
$0<\alpha<1$ is
\begin{gather*}
D^{\alpha}_{a+}f(t)=\frac{1}{\Gamma(1-\alpha)}\frac{d}{dt}\int^{t}_a
\frac{f(\tau)\,d\tau}{(t-\tau)^{\alpha}}.
\end{gather*}
The right-side fractional derivative is
\begin{gather*}
D^{\alpha}_{b-}f(t)=
\frac{1}{\Gamma(1-\alpha)}\left(-\frac{d}{dt}\right)\int^{b}_t\frac{f(\tau)\,d\tau}{(\tau-t)^{\alpha}}.
\end{gather*}
The left-side Riemann--Liouville fractional derivative of the arbitrary
order $\alpha\in\mathbb{R}$ is
\begin{gather*}
D^{\alpha}_{a+}f(t)=\frac{1}{\Gamma(n-\alpha)}\left(\frac{d}{dt}\right)^{n}\int^{t}_a
\frac{f(\tau)\,d\tau}{(t-\tau)^{1+\alpha-n}}.
\end{gather*}
The right-side fractional derivative is
\begin{gather*}
D^{\alpha}_{b-}f(t)=
\frac{1}{\Gamma(n-\alpha)}\left(-\frac{d}{dt}\right)^{n}
\int^{b}_t\frac{f(\tau)\,d\tau}{(\tau-t)^{1+\alpha-n}},
\end{gather*}
where $n=[\alpha]+1$.

The corresponding \textit{regularized fractional derivative} of
the function $f(x)$, or its Caputo fractional derivatives (in
honour of M.~Caputo~\cite{Capuro}),
are def\/ined as follows
\begin{gather*}
{}^cD^{\alpha}_{a+}f(t)=\frac{1}{\Gamma(1-\alpha)}
\left[\frac{d}{d\tau}\int^{t}_a\frac{f(\tau)\,d\tau}{(t-\tau)^{\alpha}}-\frac{f(a)}{(t-a)^{\alpha}}\right],
\end{gather*}
the right-side derivative being
\begin{gather*}
{}^cD^{\alpha}_{b-}f(t)=\frac{1}{\Gamma(1-\alpha)}
\left[\left(-\frac{d}{d\tau}\right)
\int^{b}_t\frac{f(\tau)\,d\tau}{(\tau-t)^{\alpha}}-\frac{f(b)}{(b-t)^{\alpha}}\right].
\end{gather*}

The Caputo fractional derivative of the arbitrary order $\alpha
\in\mathbb{R}$ is
\begin{gather}
\label{fracalc:CL}
{}^cD^{\alpha}_{a+}f(t)=\frac{1}{\Gamma(n-\alpha)}
\left[\left(\frac{d}{d\tau}\right)^{n}
\int^{t}_a\frac{f(\tau)\,d\tau}{(t-\tau)^{1+\alpha-n}}-\frac{f(a)}{(t-a)^{\alpha}}\right],
\end{gather}
and the right-side derivative is
\begin{gather}
\label{fracalc:CR}
{}^cD^{\alpha}_{b-}f(t)=\frac{1}{\Gamma(n-\alpha)}
\left[\left(-\frac{d}{d\tau}\right)^{n}
\int^{b}_t\frac{f(\tau)\,d\tau}{(\tau-t)^{1+\alpha-n}}-\frac{f(b)}{(b-t)^{\alpha}}\right].
\end{gather}

If $f(x)$ is absolutely continuous on $[a,T]$ or corresponding
$[T,b]$, then the left Caputo fractional derivative is
\begin{gather}
\label{fracalc:03}
{}^cD^{\alpha}_{a+}f(t)=\frac{1}{\Gamma(n-\alpha)}
\int^{t}_a\frac{d\tau}{(t-\tau)^{1+\alpha-n}}\left(\frac{d}{d\tau}\right)^n f(\tau),
\end{gather}
the right-side Caputo fractional derivative being
\begin{gather}
\label{fracalc:04}
{}^cD^{\alpha}_{b-}f(t)=\frac{1}{\Gamma(n-\alpha)}
\int^{b}_t\frac{d\tau}{(\tau-t)^{1+\alpha-n}}\left(-\frac{d}{d\tau}\right)^n f(\tau),
\end{gather}
where~$\alpha$ represents the order of the derivative,
$n-1<\alpha<n$.

The forms~\eqref{fracalc:03} and~\eqref{fracalc:04} are often used
in physical literature.
However, the forms~\eqref{fracalc:CL} and
\eqref{fracalc:CR} are more useful, since they may be applied to a
wider class of functions.

The Caputo fractional derivative or integral in the Weyl sense is
the corresponding Caputo operator for the absolutely continuous
function $f(x)$ def\/ined on the whole real axes $\mathbb{R}$.

Some properties of fractional derivatives and~integrals are listed
below
\begin{gather}
D^{-\alpha}_{a+}=I^{\alpha}_{a+}f(t)\qquad
\big(D^{-\alpha}_{b-}=I^{\alpha}_{b-}f(t)\big),\qquad\alpha>0,
\nonumber
\\
I^{\alpha}_{a+}f(t)=
\frac{1}{\Gamma(\alpha)}\int^{t}_a\frac{f(\tau)\,d\tau}{(t-\tau)^{\alpha-1}},\qquad t>a,
\nonumber
\\
I^{\alpha}_{b-}f(t)=\frac{1}{\Gamma(\alpha)}\int^{b}_t\frac{f(\tau)\,d\tau}{(\tau-t)^{\alpha-1}},\qquad t<b.
\nonumber
\\
\label{fracalc:07}
D^{\alpha}_{a+}f(t)=I^{-\alpha}_{a+}f(t),
\\
D^{\alpha}_{a+}D^{\beta}_{a+}f(t)=D^{\beta}_{a+}D^{\alpha}_{a+}f(t)=
D^{\alpha+\beta}_{a+}f(t),
\nonumber
\\
I^{\alpha}_{a+}I^{\beta}_{a+}f(t)=I^{\beta}_{a+}I^{\alpha}_{a+}f(t)=
I^{\alpha+\beta}_{a+}f(t),
\nonumber
\\
f(t)=\sum^{n-1}_{j=0}\frac{D^{\alpha+j}_{a+}f(0)}{\Gamma(1+\alpha+j)}t^{\alpha+j}+R_n(t),\qquad n=
[\textrm{Re}\,\alpha]+1,
\nonumber
\end{gather}
where $R_n(t)=I^{\alpha+j}_{a+}D^{\alpha+j}_{a+}f(t)$.

The derivatives of some functions are:
\begin{gather*}
{}_{-\infty}D^{\alpha}_{t+}\sin{\lambda
t}=\lambda^{\alpha}\sin{\left(\lambda t+\frac{\pi\alpha}{2}\right)},
\qquad
{}_{-\infty}D^{\alpha}_{t+}\cos{\lambda
t}=\lambda^{\alpha}\cos{\left(\lambda t+\frac{\pi\alpha}{2}\right)},
\end{gather*}
where $\lambda>0,\alpha>-1$; when $\alpha\leq-1$, we have to
use the property~\eqref{fracalc:07};
\begin{gather*}
{}_{-\infty}D^{\alpha}_{t+}\textrm{e}^{\lambda t+\mu}
=\lambda^{\alpha}\textrm{e}^{\lambda t+\mu},
\qquad
\operatorname{Re} \lambda>0.
\end{gather*}

The Riesz fractional derivative and~integral of order~$\alpha$ are
def\/ined by the Fourier transformation
\begin{gather*}
\textrm{D}^{\alpha}_{x}f(x)=\mathcal{F}^{-1}\left(|k|^{\alpha}(\mathcal{F}f)(k)\right),
\qquad
\textrm{I}^{\alpha}_{x}f(x)=\mathcal{F}^{-1}\left(|k|^{-\alpha}(\mathcal{F}f)(k)\right).
\end{gather*}

The Riesz fractional integral could be presented as a~convolution
\begin{gather*}
\textrm{I}^{\alpha}_{x}f(x)=\int_{\mathbb{R}^n}K_{\alpha}(x-\xi)f(\xi)\,d\xi,\qquad
\alpha>0,
\end{gather*}
where $K_{\alpha}(x)$ is the Riesz kernel
\begin{gather*}
K_{\alpha}(x)=
\begin{cases}
\gamma^{-1}_{n}(\alpha)|x|^{\alpha-n}, &\alpha\not = n+2k,
\\
-\gamma^{-1}_{n}(\alpha)|x|^{\alpha-n}\ln{|x|},\qquad &\alpha=n+2k,
\end{cases}
\end{gather*}
where $k\in\mathbb{N}$ and~coef\/f\/icients
$\gamma^{-1}_{n}(\alpha)$ are
\begin{gather*}
\gamma^{-1}_{n}(\alpha)=
\begin{cases}
2^{\alpha}\pi^{n/2}\Gamma(\frac{\alpha}{2})/\Gamma\left(\frac{n-\alpha}{2}\right), &\alpha\not = n+2k,
\\[0.9mm]
(-1)^{(n-\alpha)/2}2^{\alpha-1}\pi^{n/2}\Gamma(\frac{\alpha}{2})\Gamma\left(1+\frac{\alpha-n}{2}\right),\qquad &\alpha=n+2k.
\end{cases}
\end{gather*}

Note that the correlation of the Riesz and~Caputo fractional derivatives in the Weyl~sense
\begin{gather*}
\textrm{D}^{\alpha}_{x}f(x)=\frac{1}{2\cos{(\alpha\,\pi/2)}}
\left({}^cD^{\alpha}_{+}f(x)+{}^cD^{\alpha}_{-}f(x)\right).
\end{gather*}

Two special functions of those often used in fractional calculus are as follows:

the Mittag--Lef\/f\/ler function
\begin{gather*}
\qquad E_{\alpha,\beta}(z)=
\sum^{\infty}_{n=0}\frac{z^{n}}{\Gamma(\alpha n+\beta)},
\end{gather*}

the generalized exponential function
\begin{gather*}
\qquad E^z_{\alpha}=\sum^{\infty}_{n=0}\frac{z^{n+\alpha}}{\Gamma(1+\alpha+n)}.
\end{gather*}

The basic aspects of fractional calculus and~its applications could be
found in~\cite{SamkoKilbasMarichev1993,Podlubny1999,Gerasimov,Nigmatullin1984,Miskinis2003,
MillerRoss1993,KlafterShlesingerZumofen1996,CarpinetryMainardy1997,Bardouetal2002,Hilfer2002,Zaslavsky2002}.

\section{Supersymmetry}
\label{appendixC}

{\sloppy Let us consider an important example.
Two large inf\/inite classes of
groups, ${\rm Osp}(N|M)$ and ${\rm SU}(N|M)$, are used in applications.
The orthogonal group ${\rm O}(N)$ preserves the inva\-riant~$x_ix^i$, and
the group ${\rm Sp}(M)$ retains the invariant $\theta_mC_{mn}\theta^n$, where the
$C_{mn}$ matrices are real antisymmetric matrices and
$\theta_i$ are Grassmann-valued.
The orthosymplectic group is
now def\/ined as the group that preserves the sum
\begin{gather*}
{\rm Osp}{(N|M)}:\
x_ix^i+\theta_mC_{mn}\theta^n=\operatorname{inv}.
\end{gather*}
Note that the orthosymplectic group contains the product
\begin{gather}
{\rm Osp}{(N|M)}=\left(
\begin{matrix}
{\rm O}(N)&A
\\
B&{\rm Sp}(M)
\end{matrix}\right).
\end{gather}

}

Let us write the generators of ${\rm Osp}(1|4)$ as
$M_A=(p_{\mu},M_{\mu,\nu},Q_{\alpha})$ which have the
commutation relations
\begin{gather}
[M_A,M_B]_{\pm}=f^C_{AB}M_C.
\end{gather}
In the explicit form, the commutators involving the supersymmetry
generator are
\begin{gather*}
\{Q_{\alpha},Q_{\beta}\}=2(\gamma^{\mu}C)_{\alpha\beta}P_{\mu},
\qquad
\left[Q_{\alpha},P_{\mu}\right]=0,
\qquad
\left[Q_{\alpha},M_{\mu,\nu}\right]=(\sigma_{\mu,\nu})^{\beta}_{\alpha}Q_{\beta}.
\end{gather*}
Let us def\/ine superspace as the space created by the pair
$
x_{\mu}\,\theta_{\alpha}$,
where $\theta_{\alpha}$ is a~Grassmann number.
Let us def\/ine
the supersymmetry generator as
\begin{gather*}
Q_{\alpha}=
\frac{\partial}{\partial\bar{\theta}^{\alpha}}-i(\gamma^{\mu}\theta)_{\alpha}\partial_{\mu}.
\end{gather*}
The anticommutator between two such generators yields a
displacement
\begin{gather*}
\{Q_{\alpha},Q_{\beta}\}=-2(\gamma^{\mu}C)_{\alpha\beta}i\partial_{\mu}.
\end{gather*}
Note that $\bar{\varepsilon}Q$ makes the following
transformations of the superspace
\begin{gather*}
x_{\mu}\rightarrow
x_{\mu}-i\bar{\varepsilon}\gamma_{\mu}\theta,
\qquad
\theta_{\alpha}\rightarrow\theta_{\alpha}+
\varepsilon_{\alpha}.
\end{gather*}
We can construct the operator
\begin{gather*}
D_{\alpha}=
\frac{\partial}{\partial\bar{\theta}^{\alpha}}+i(\gamma^{\mu}\partial_{\mu}\theta)_{\alpha}.
\end{gather*}
This operator $D_{\alpha}$ anticommutes with the supersymmetry
generator,
$\{Q_{\alpha},D_{\beta}\}=0$.
This relation is very important because it allows us to place
restrictions on the representations of supersymmetry without
destroying the symmetry.
This permits us to extract the
irreducible representations from the reducible ones.

\subsection*{Acknowledgements}

The author would like to express his gratitude to Professors B.A.~Dubrovin, M.~Pavlov
and L.~Ala\-niya for the invitation and~kind hospitality during the Conference ``Geometrical Methods in Mathematical Physics''
(Moscow State University, December 12--17, 2011).

\pdfbookmark[1]{References}{ref}
\LastPageEnding

\end{document}